\documentclass[10pt,hidelinks,twocolumn]{amsart}
\usepackage{fullpage}
\usepackage{color}

\usepackage[english]{babel}
\usepackage{hyperref}
\input xy
\xyoption{all}

\usepackage{color, graphicx}
\usepackage{amsmath}
\usepackage{amssymb}
\usepackage{wrapfig}

\usepackage{calrsfs}

\theoremstyle{plain}

\newtheorem{theorem}{Theorem}

\theoremstyle{definition}
\newtheorem{definition}[theorem]{\bf Definition}

\theoremstyle{remark}

\newcommand{\HIDE}[1]{}

\newcommand{\N}{{\mathbb N}}
\newcommand{\Q}{{\mathbb Q}}

\newcommand{\Z}{{\mathbb Z}}

\newcommand{\be}{\begin{enumerate}}
\newcommand{\ee}{\end{enumerate}}

\newcommand\cyr{%
\renewcommand\rmdefault{wncyr}%
\renewcommand\sfdefault{wncyss}%
\renewcommand\encodingdefault{OT2}%

\normalfont 
\selectfont}
\DeclareTextFontCommand{\textcyr}{\cyr}

\title{\Large In memory of Martin Davis}
\author{Wesley Calvert}
\author{Valentina Harizanov}
\author{Eugenio G. Omodeo}
\author{Alberto Policriti}
\author{Alexandra Shlapentokh}


\thanks{}

\graphicspath{{./images/}}

\begin{document}
\maketitle












\begin{sloppypar}
\noindent In 1950, Martin David Davis found the culture in his graduate program at Princeton deeply alienating, and he was ready to be out of there.  So he solved an open problem of Kleene (establishing the backbone for a major branch of modern computability theory), wrote down his initial steps toward the eventually successful solution of a Hilbert problem, and graduated.

\begin{figure}[!htbp]
    \centering
    \includegraphics[width=0.4\textwidth]{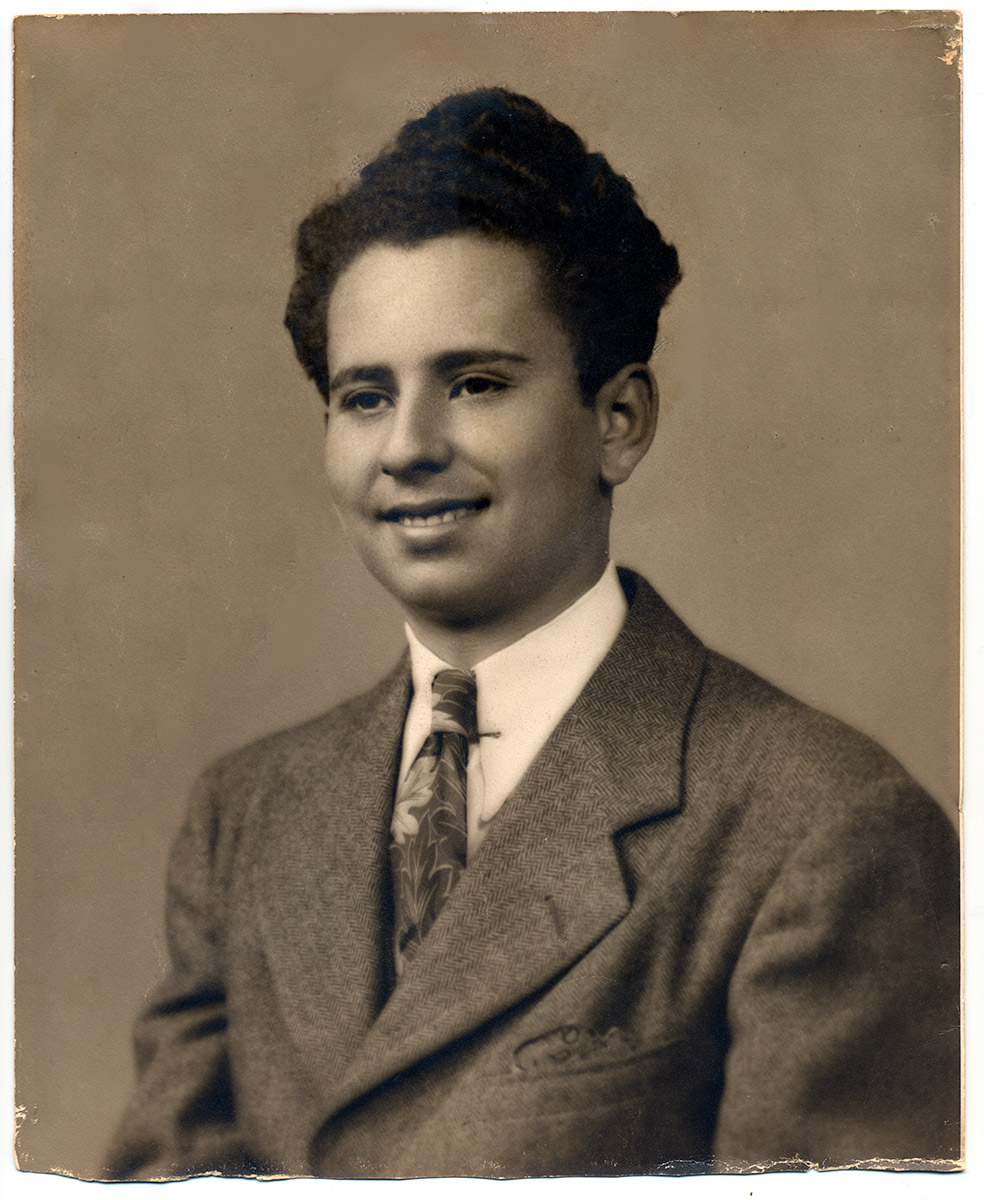}
    \caption{Martin in the late 1940s.  (Courtesy of his son Harold Davis)}
    \label{fig:youngMartin}
\end{figure}

Davis's solution to the problem of Kleene became the hyperarithmetical hierarchy, which we will explain in Section \ref{hypsection}.  His work on Hilbert's Tenth Problem included what would later be called his ``daring hypothesis" \cite{Matiyasevich93}: a conjecture, later verified, that the computably enumerable sets were exactly the Diophantine sets, as we will descuss in Section \ref{H10sec}.

In spite of the economic restrictions which his family, as immigrants who arrived shortly before the Great Depression, had to face, he received a high-quality education. He arrived at  City College of New York as a freshman in 1944 and soon became interested in the foundations of real analysis and in logic. He approached Post, who introduced him to the writings by Church and Kleene on algorithmic unsolvability and to Hilbert's Tenth Problem, which would soon become Davis's ``lifelong obsession''.  When Davis had to choose where to undertake graduate studies, Post advised him to go to Princeton, where, as Davis later expressed it, the ``culture clash" between his Jewish working class background and the ``genteel Princeton atmosphere" made him eager to conclude quickly: in fact, he got his Ph.D.\ in just 2 years, under the guidance of Church, in 1950. 

\begin{figure}[!htbp]
    \centering
    \includegraphics[width=0.4\textwidth]{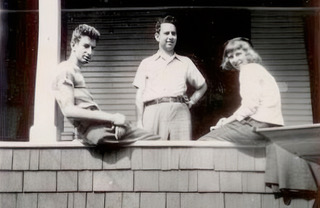}
    \caption{Martin, middle, with his friends Jacob T.\ and Judith Schwartz. (Courtesy of Diana Robinson Schwartz)}
    \label{fig:MartinWJackwJudith}
\end{figure}


Davis's first position was at the University of Illinois at Urbana-Champaign, but ``the Korean war and the hot breath of the draft'' led him to leave that job for the Control Systems Laboratory.  He later moved to the Institute for Advanced Study, the University of California at Davis, the Ohio State University, the Rensselaer Polytechnic Institute, Yeshiva University, and New York University. Certain summer projects funded by military and civilian research agencies enabled him to make crucial achievements (``It was in the summer of 1959 that Hilary and I really hit the
jackpot,'' he says, to describe the original, raw discovery of what would become known as the celebrated Davis--Putnam--Robinson theorem).  Over the course of his career, Davis supervised a total of 25 Ph.D.\ students, including scholars now known for work in mathematics, computer science, and philosophy.

\begin{figure}[!htbp]
    \centering
    \includegraphics[width=0.4\textwidth]{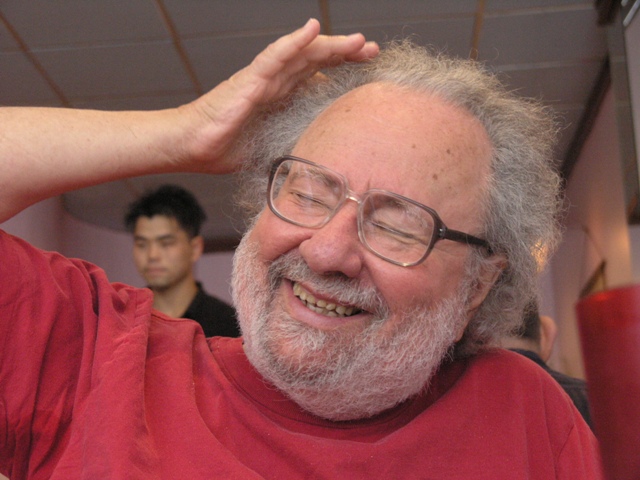}
    \caption{Martin ready to give an answer, late 1990s. (Courtesy of Domenico Cantone)}
    \label{fig:MartinReadyToAnswer}
\end{figure}

Davis's expository books have become classics and have been translated into various languages: \textit{Computability and Unsolvability}; \textit{A First Course in Functional Analysis}, \textit{Applied Nonstandard Analysis}; \textit{Computability, Complexity, and Languages: Fundamentals of Theoretical Computer Science} (with Ron Sigal and Elaine J. Weyuker); and \textit{The Universal Computer: The Road from Leibniz to Turing}. His 1993 \textit{Lecture Notes in Logic} are a true jewel.


His many honors and awards include the Steele Prize, the Chauvenet Prize (with Reuben Hersch), Fellowship of the AAAS, a Guggenheim Foundation Fellowship,
the Herbrand Award of the International Conference on Automated Deduction; and the Pioneering Achievement Award from the ACM SIG on Design Automation.

Martin Davis left two rather comprehensive autobiographic accounts \cite{Dav16,Davis20} and the long interviews \cite{AllynJackson08,Jackson2021}.  For this reason, the present note will primarily focus on his scientific achievements.

\section{Computability}

\subsection*{Computably Enumerable Sets and Universal Turing Machines.}
Modern computability theory started in 1936 with Turing's seminal paper on computable numbers with an application to the Entscheidungsproblem, a decision problem of Hilbert and Ackermann for which Turing provided a negative solution. Turing introduced what today we would call a Turing machine (he called it an a-machine), which is essentially an abstraction of a computer. Turing's and other formalisms for an intuitive concept of an effectively calculable function, developed by G\"{o}del, Kleene, Church, Post and others, had profound significance
for the emerging science of computing.


One of the main concepts in computability theory (also called recursion theory) is that of a computable function and computable relation.  
A function $f :\mathbb{N}^{n} \rightarrow \mathbb{N}$ is \emph{computable} if there is a Turing machine that on every input $a_{1} ,\ldots  ,a_{n}$ halts and outputs its value $f (a_{1} ,\ldots  ,a_{n})$. Addition and multiplication on natural numbers are computable functions (operations). A set of natural numbers is \emph{computable} if its characteristic function is computable. For example, the set of prime numbers is computable. All finite sets are computable. Clearly, the complement of a computable set is computable.  Computable $n$-ary relations are defined similarly. It can be shown that there is a computable bijection $h :\mathbb{N}^{2} \rightarrow \mathbb{N}$, thus allowing algorithmic coding of pairs and, more generally, finite tuples of natural numbers by natural numbers. Decidable problems are encoded
by computable relations. 

A function $f :D \rightarrow \mathbb{N}$, where $D\subseteq\mathbb{N}^n$, is \emph{partial computable} if there is a Turing machine that on every input in the domain $D$ of $f$ halts and outputs its value, while on every input in $\mathbb{N}^{n}$ that is not in the domain of $f$ it does not halt, thus computing forever. Clearly, computable functions are partial computable functions that are total. Partial computable functions coincide with partial recursive functions defined by Kleene, starting with some basic functions and applying the operations of composition, primitive recursion, and unbounded search.

Since each Turing machine is a finite list of instructions, Turing machines can be algorithmically enumerated without repetitions as: \medskip 

$M_{0} ,M_{1} ,M_{2} ,\ldots.
$\medskip

\noindent A Turing machine on a given input may halt and output its value, or it may compute forever. For each Turing machine $M_{e}$, we denote the $n$-ary partial function it computes by $\varphi _{e}^{(n)}$ and its domain by $W_{e}^{(n)}$. Hence \medskip

$\varphi _{0}^{(n)} ,\varphi _{1}^{(n)} ,\varphi _{2}^{(n)} ,\ldots $ \medskip

\noindent is a computable enumeration of all $n$-ary partial computable functions. For $n =1 ,$ we omit the superscript. Moreover, there is a binary partial computable function $\psi $ such that $\psi (e ,x) =\varphi _{e}(x)$. \medskip

The above enumeration gives rise to a universal Turing machine, which can simulate any Turing machine on any input and leads to the idea of stored-program computer. Davis has written on the universal Turing machine in a number of papers starting in 1956. His lecture in 2012 titled ``Universality is ubiquitous'' is available at \url{https://www.youtube.com/watch?v=ZVTgtODX0Nc}

In a written version published in \cite{Davis13b}, Davis wrote: ``Turing's concept of `universal machine' will be discussed as an abstraction,
as embodied in physical devices, as present in nature, and in connection with the artificial intelligence project''.

In addition to his large body of expository work on universal Turing machines, Davis was also an early technical contributor to the subject.  Turing had constructed a universal machine, but had not dealt with universal Turing machines as a class of objects.  John McCarthy and Claude Shannon posed the problem of giving a definition of universal Turing machines, which would deal, for instance, with the simplicity of the encoding by which the universal machine simulates arbitrary machines.  Davis solved this problem in \cite{Davis56}, as we explain next.  

\begin{definition} Let $S$ be a set.  We say that $S$ is \emph{computably enumerable} (also called \emph{recursively enumerable}) if and only if $S$ is empty or the range of a computable function.
    \end{definition}

It is not hard to see that a set is computably enumerable if and only if it is the domain $W_{e}$ of some partial computable (equivalently, partial recursive) function $\varphi _{e}$.  If $W_{e}$ is nonempty, then $W_{e}$ can be computably enumerated by the procedure that simultaneously
runs \medskip

\[M_{e} (0) ,M_{e} (1) ,\ldots  ,M_{e} (k) ,\ldots \] 

 \noindent 
 and enumerates those $k$ for which $M_{e} (k)$ halts, as soon as the halting occurs. Here, simultaneously means that at each step we add a new input and also
run all activated inputs for an additional computational step.
The converse, that every computably enumerable set is some $W_{e}$, is also true.

Since we can computably enumerate the Turing machines, we can also computably enumerate the computably enumerable sets by 
\[W_{0} ,W_{1} ,W_{2} ,\ldots.\] 

Clearly, every computable set is computably enumerable since a decision algorithm can be transformed into an enumeration algorithm. Computable sets are exactly computably enumerable sets that also have computably enumerable complements.

\begin{definition} \cite{Davis56}
    \begin{enumerate}
    \item We say that a set $S$ is \emph{c.e.\ complete} if and only if it is computably enumerable and for any computably enumerable set $W$, there is a computable function $\sigma:\mathbb{N} \to \mathbb{N}$ with $x \in W$ if and only if $\sigma(x) \in S$.
    \item We define $\delta_M$ to be the set of all initial configurations of $M$ from which $M$ will eventually halt. 
    \item We say that $M$ is \emph{universal} if and only if $\delta_M$ is a computably enumerable set that is complete.
\end{enumerate}
\end{definition}

Davis proved that a machine which is universal in this sense does, indeed, simulate all Turing machines.  However, having universal machines defined as a class, rather than simply observed as a phenomenon, opened the door to lines of thinking that involve quantification over all universal Turing machines, such as the foundational work of Solomonoff \cite{Solomonoff1964} and Kolmogorov \cite{Kolmogorov1965} on information theory (Kolmogorov complexity), and the theory of algorithmic randomness arising, in part, from it \cite{DowneyHirschfeldt2010}.


The diagonal \emph{halting set}
$H$ consists of all inputs $e$ on which the Turing machine with index $e$ halts. That is, \medskip  
\begin{center}
$H =\{e$$ :M_{e}$ halts on input $e\}$.
\end{center}\par

\noindent It can be shown that the set $H$ is computably enumerable. It is not computable since its complement is not computably enumerable. If the complement $\overline{H}$ were computably enumerable, then for some $e_{0}$, we would have $\overline{H} =W_{e_{0}}$. Then \medskip 
\begin{center}
$e_{0} \in H \Leftrightarrow e_{0} \in W_{e_{0}} \Leftrightarrow e_{0} \in \overline{H}\text{,}$ \end{center}\par

\noindent which is a contradiction.

It can also be shown that a set $A$ is computably enumerable\ if and only there is a computable binary relation $R$ such that for every $a$,  \medskip

\begin{center}
$a \in A \Leftrightarrow ( \exists x)R (a ,x)\text{.}$ \medskip
\end{center}\par

We can relativize all of these notions using Turing's notion of an \emph{oracle} machine.  A machine with an oracle for a set $S$ is a Turing machine which carries out its computation with the additional resource of read-only access to the characteristic function of the set $S$.  In this way, even if a set $U$ is not computable, it may be computable with an oracle for another set --- for instance, its complement.  The halting set relative to $S$ (also called the \emph{jump} of $S$ and denoted $S'$) is defined exactly as before, but replacing the machines with machines with oracle $S$.

\subsection*{The Hyperarithmetical Hierarchy}\label{hypsection}

In his 1950 Ph.D.\ thesis at Princeton, Martin addressed a problem posed by Kleene.  It was already known that every formula of classical predicate logic is equivalent to a formula consisting of a block of quantifiers (``for all'' and ``there exists''), followed by a quantifier-free formula.  Kleene noted that the optimal form of such an equivalent formula corresponded to the degree of unsolvability of satisfying that formula.  For instance, Turing's halting problem is equivalent to the problem of satisfying a particular sentence with a single existential quantifier, but not to the satisfaction of any quantifier-free formula.

This gives rise to a hierarchy of formulas --- and, equivalently, of decision problems, according to the number of quantifiers, and whether those quantifiers are ``for all'' or ``there exists.''  Since there is no computational difference between determining the existence of a single element and determining the existence of a finite tuple of elements, we consider only alternations of quantifiers.

\begin{definition}[The Arithmetical Hierarchy] Let $S \subseteq \mathbb{N}^m$.
\begin{enumerate}
    \item We say that $S$ is $\Sigma^0_1$ if and only if there is a computable set $T \subseteq \mathbb{N}^{m+1}$ such that $\bar{a} \in S$ if and only if $\exists x \left[ (x,\bar{a}) \in T\right]$.
    \item We say that $S$ is $\Pi^0_1$ if and only if there is a computable set $T \subseteq \mathbb{N}^{m+1}$ such that $\bar{a} \in S$ if and only if $\forall x \left[ (x,\bar{a}) \in T\right]$.
    \item We say that $S$ is $\Sigma^0_{n+1}$ if and only if there is a $\Pi^0_{n}$ set $T \subseteq \mathbb{N}^{m+1}$ such that $\bar{a} \in S$ if and only if $\exists x \left[ (x,\bar{a}) \in T\right]$.
    \item We say that $S$ is $\Pi^0_{n+1}$ if and only if there is a $\Sigma^0_{n}$ set $T \subseteq \mathbb{N}^{m+1}$ such that $\bar{a} \in S$ if and only if $\forall x \left[(x,\bar{a}) \in T\right]$.
\end{enumerate}
\end{definition}

The $\Sigma^0_1$ sets are exactly the computably enumerable sets.  The $\Sigma^0_{n+1}$ sets are exactly those computably enumerable relative to the $n$-times iterated jump of the empty set.  There are certainly sets of natural numbers that are not $\Sigma^0_n$ or $\Pi^0_n$ for any $n$, and Kleene asked whether the hierarchy could be continued to transfinite levels.

Davis \cite{Davis50a} carried out this generalization by iterating the jump starting at an arbitrary set, and also describing a uniform join of infinitely many jumps.  In the following definition, $\omega$ denotes the least transfinite ordinal, the order type of the natural numbers.

\begin{definition} Let $K_0 = \emptyset$.
\begin{enumerate}
    \item Let $K_{\alpha+1} = K_\alpha'$.
    \item Let $K_{\omega n}$ be the set defined by $2^{x_1}3^{x_2} \in K_{\omega n}$ if and only if $x_1 \in K_{\omega (n-1) +x_2}$.
\end{enumerate}
\end{definition}

Davis proved that this hierarchy is proper and that it extends the finite-level Kleene hierarchy to all ordinals less than $\omega^2$.  After Davis's work, the major pieces still missing were extension to larger ordinals and the fact that the choice of representative sets at limit levels is not unique.  It was five years later that Spector showed that, up to Turing degree, the definition is robust \cite{Spector1955}.  This transfinite hierarchy is now at the core of modern computable mathematics \cite{ashknight}.

\subsection*{Hypercomputation, Neural Networks, and Unconventional Computation}

A major effort of the later part of Davis's career was devoted to defining the bounds of computation.  In response to a community of scholars who applied relativistic and quantum theories to propose computing devices more powerful than a Turing machine (an idea they called ``hypercomputation''), Davis responded with unbridled skepticism.

Davis's key paper \cite{Davis04} on this subject  considered a key proposal for a hypercomputer, a certain kind of neural network.  The proposal described a certain model with parameters.  If those parameters range over rational numbers, the machine could determine membership in the computable sets, as expected; however, if the parameters were allowed to range over arbitrary reals, it could determine membership in arbitrary subsets of $\mathbb{N}$.  Davis pointed out that the ability of the model to determine membership in arbitary sets followed immediately from Turing's observation that any countable set was computable relative to some oracle, still via a Turing machine.  In this way, Davis showed that neural networks reflect computation within the bounds of the Church--Turing Thesis --- that is, equivalent to a computation that can be done by a Turing machine.

Many of the hypercomputation models rely on some access to certain full-precision real numbers, to which, Davis pointed out, no scientific observation could give us access.

In a paper entitled, ``Why there is no such discipline as hypercomputation'' (published as an introduction to a special issue of \emph{Applied Mathematics and Computation} devoted to papers on exactly this discipline), Davis argued that if there were a hypercomputer, we would be unable to verify its performance, since we could only see finitely many outputs.  Moreover, any real computer is subject not merely to the limitation of a Turing machine --- that is, the limitation that only a finite time and a finite span of memory can be used --- but to a much more strict limitation of a constant bound on these quantities, depending only on the machine, and not on the algorithm or the data.  While investigation of algorithms in the Turing machine context has important meaning, both in theory and in practice, they must finally be executed on finite state automata, a much more limiting device.

\subsection*{On a personal note VH} The first time I\  met Martin Davis at a conference, I was very impressed by his kindness, modesty, sense of humor, and friendliness. Later, I invited him and Virginia to visit me at George Washington University and give a math colloquium talk. It was in November 2007 and his lecture was titled ``Unsolvability and undecidability in the Diophantine realm,'' an ordinary title compared to his 2020 MSRI talk ``Here there be monsters.'' During the GW talk to a packed room, Davis covered many years of work, progress, and struggles on Hilbert's Tenth Problem. When a person in the audience asked him how he managed to persist for so many years working on one problem, he replied that his obsession with this problem ``was a disease.'' \medskip

Because of the advances of computer technology, Martin Davis's 2020 talk, ``Here there be monsters,'' was possible when MSRI\ semester program on Decidability, Definability and Computability in Number Theory had to be moved online. His talk was an opening one for the program and is posted on \url{https://www.msri.org/seminars/25120}.

\begin{figure}[htb!]
\centering
    \includegraphics[width=0.4\textwidth]{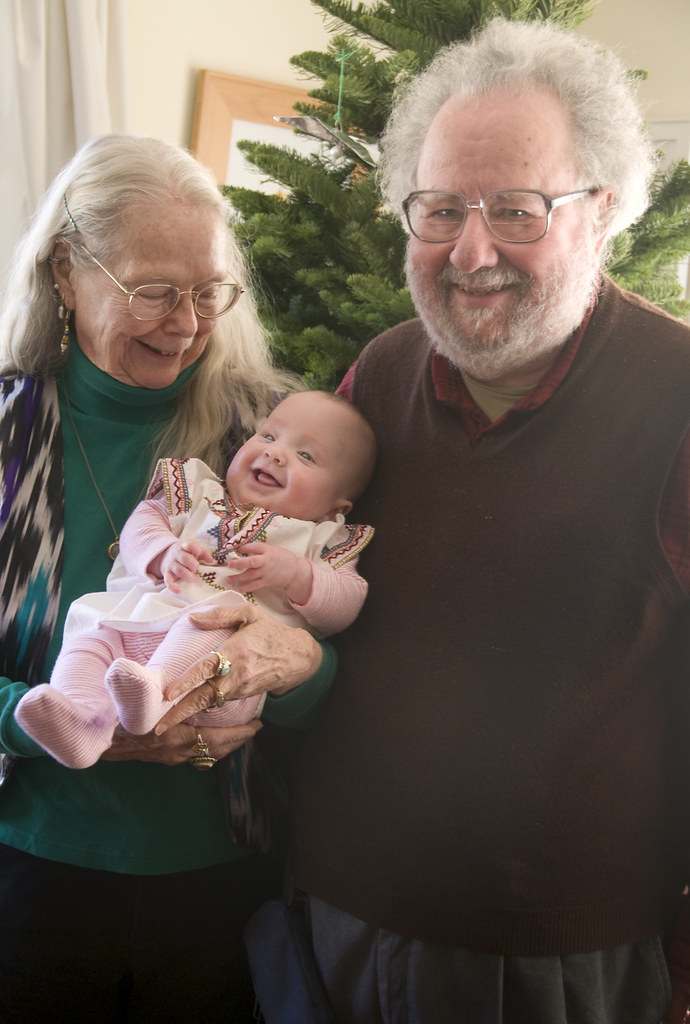}
\caption{\label{fig:Martin_w_Virginia_grandparents}Martin and his wife Virginia as grandparents, holding Katie Rose. (Courtesy of Harold Davis)}
\end{figure}

\section{Hilbert's Tenth Problem}\label{H10sec}

\subsection*{A brief history} Martin Davis was one of the four people who collectively solved Hilbert's Tenth Problem.  The other three were Hilary Putnam, Julia Robinson and Yuri Matiyasevich. 

The history of Hilbert's Tenth Problem starts in 1900 when, during an International Congress of Mathematicians in Paris, David Hilbert presented a list of 23 problems that had a great influence on mathematics in the twentieth century and continue to influence the subject in the twenty-first. The tenth problem on the list asked to devise a process that determines whether any given Diophantine (polynomial) equation with integer coefficients
has a solution in the integers.

If we are to rephrase Hilbert's question in modern terms, we could say that he asked for an algorithm (or a computer program) taking as its input coefficients of a polynomial equation in several variables over $\Z$ and generating a ``yes'' or ``no'' answer  to the question about existence of the roots of this polynomial over $\Z$.  

At the time Hilbert formulated his question a formal notion of an algorithm, let alone computer program, did not yet exist.  He asked for a process terminating in a finite number of steps, and this was later interpreted to mean an algorithm.  The theorem proved by Davis, Putnam, Robinson and Matiyasevich showed that such an algorithm does not exist.

Lagrange's four-squares
theorem from 1772 establishes that every natural number can be expressed as the sum of squares of four integers. Hence, the algorithmic solvability of a Diophantine equation in the integers is equivalent to the algorithmic solvability of a Diophantine equation in the natural numbers.

The first step towards the solution of the problem was made by Davis in 1949 when he showed that any computably enumerable of natural numbers  has the following form:
\[
\left\{a\mid \exists y\forall k\leqslant y\exists x_{1},\ldots ,x_{n}:p\left(a,k,y,x_{1},\ldots ,x_{n}\right)=0\right\},
\]
where $p(\ldots )$ is a polynomial with coefficients in $\Z$ and all variables range over $\Z$.

It is not hard to see that a set of natural numbers defined using existential quantifiers and a polynomial equation is computably enumerable.  More precisely, let $p(t,x_1,\ldots,x_n)$ be a polynomial in $n+1$ variables and consider the following set $S$ of natural numbers:
\[
\{a \in \N|\exists x_1 \in \N \,\ldots \, \exists x_n \in \N: p(a,x_1,\ldots, x_n)=0\}.
\]
We can enumerate all $(n+1)$-tuples of natural numbers and plug them into the polynomial $p$.  Each time the result is 0, we add the first coordinate of the $(n+1)$-tuple to $S$ eventually listing every element of $S$. The polynomial $p$ is called a {\it Diophantine definition} of $S$ and sets defined using existential quantifiers and polynomial equations are called {\it Diophantine sets}.

In 1953, Davis \cite{Davis53} established that the collection of Diophantine sets is closed under unions and intersections, but not under complements. 

Observe that Davis's formula defining all computably enumerable sets is very similar to the characterization of Diophantine sets (which clearly defines at least some computably enumerable sets).  In 1949, Davis conjectured that every computably enumerable set is definable by an existential polynomial formula.  It took 20 years for this conjecture to be proven. 

\begin{figure}[htb!]
\centering
    \includegraphics[width=0.44\textwidth]{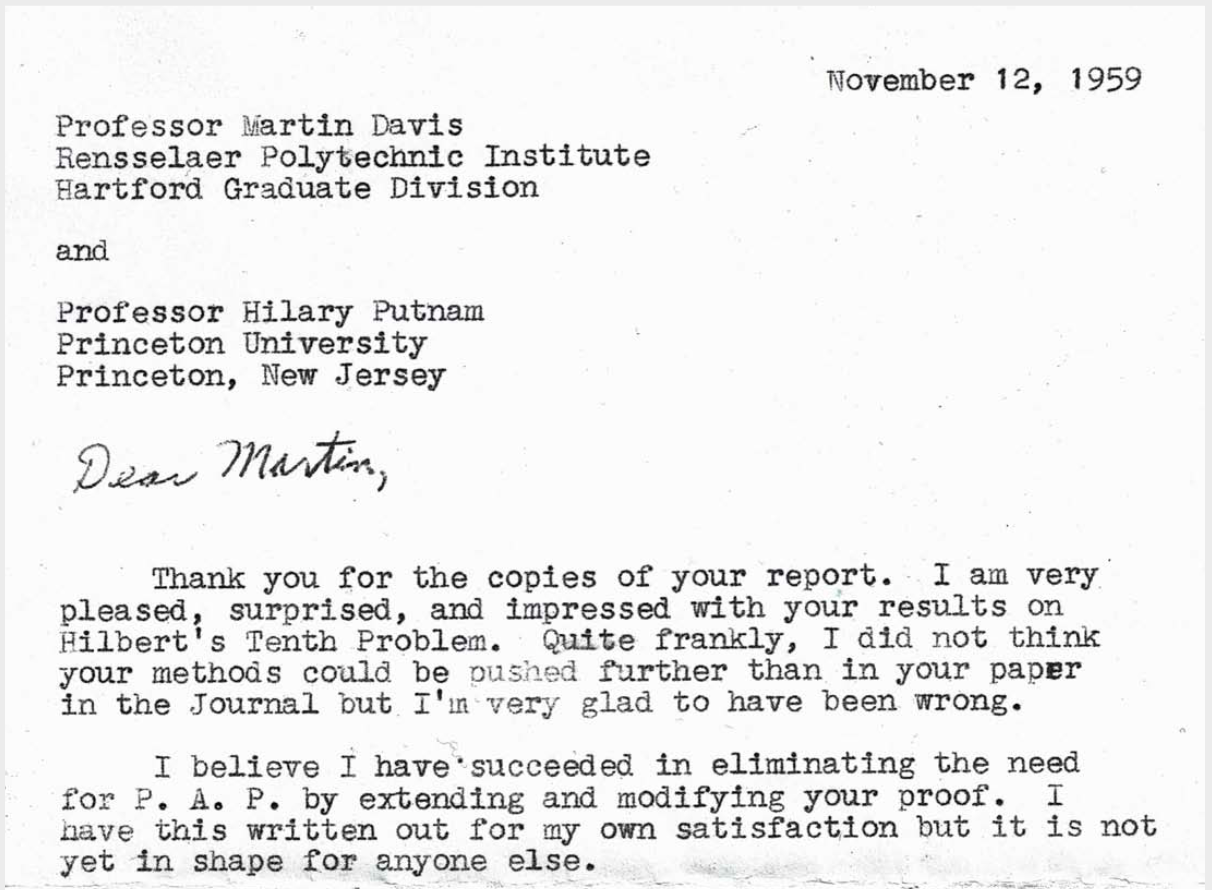}
\caption{\label{fig:Robinson_eliminates_PAP}Julia Robinson announces the possibility of getting rid of the hypothesis made by Davis and Putnam. (Courtesy of Harold Davis)}
\end{figure}

Davis and Putnam proved that, under an additional hypothesis, the bounded universal quantifier can be eliminated, in favor of an ``exponential Diophantine equation," that is, an equation where variables are allowed to appear in the exponents.  The additional hypothesis was at the time a conjecture, but is now a theorem, and concerns the lengths of arithmetic progressions of primes.  Robinson eliminated the need for the conjecture (Figure~\ref{fig:Robinson_eliminates_PAP} and \cite{DPR61}).  Finally in 1969, using Fibonacci numbers, Matiyasevich showed that exponential Diophantine equations can be replaced by polynomial equations.

It is not hard to see that the theorem proved by Davis, Putnam, Robinson and Matiyasevich implies a negative answer to Hilbert's question.  Indeed, assume that an algorithm requested by Hilbert exists.  Let the set $S$ defined above be a computably enumerable set that is not computable.  Then we can determine whether an $a \in \N$ is in $S$ by determining whether the polynomial $p(a,x_1,\ldots,x_n)$ has roots in $\N$.  Hence, if the algorithm for solving Diophantine equations exists, then there is an algorithm to determine membership in $S$.  This contradicts our assumption that $S$ is not computable.  Therefore, the algorithm requested by Hilbert does not exist.

In the Preface to 1982 Dover edition of his book \emph{Computability and Unsolvability}, Davis wrote: ``One of the great pleasures of my life came in February 1970, when I\ learned of the work of Yuri Matiyasevich which completed the proof of the crucial conjecture and thereby showed that Hilbert's Tenth Problem is recursively unsolvable.''

\begin{figure}[htb!]
\centering
    \includegraphics[width=0.44\textwidth]{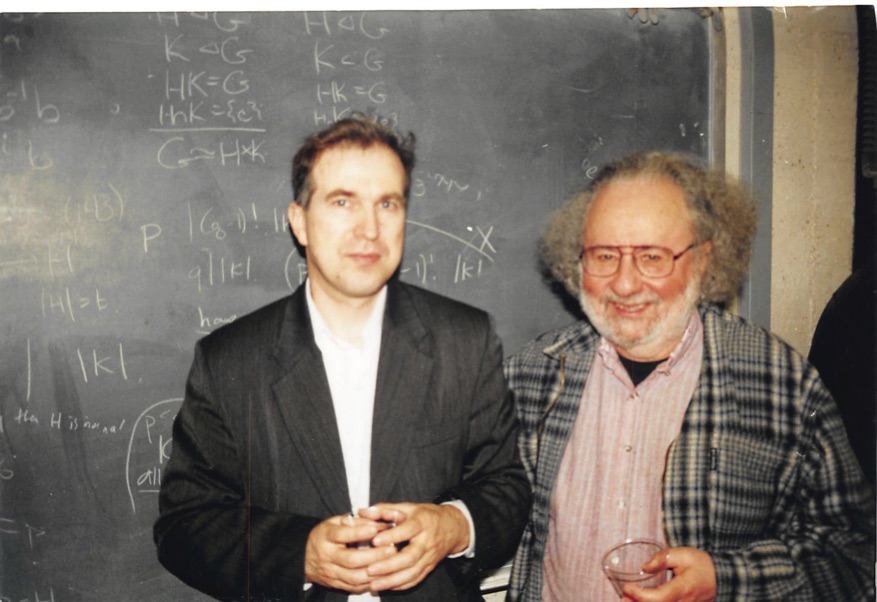}
\caption{\label{fig:MartinDavis_w_Yuri}Yuri Matiyasevich with Martin. (Courtesy of Harold Davis)}
\end{figure}

\subsection*{Other ramifications.}
The \emph{atomic diagram} of a structure $\mathcal{A}$ is the set of all atomic and negations of atomic sentences allowing additional constants for elements of the domain, which are true in $\mathcal{A}$. A structure  is \emph{computable} if the characteristic function of its atomic diagram is computable.  The standard model of arithmetic,
$\mathcal{N} =(\mathbb{N}  , + , \cdot  , S ,0)$, the natural numbers with addition, multiplication, successor function, and zero, is a computable structure. G{\"o}del established that all computable relations are definable in $\mathcal{N}$. For any computable relation there are two natural defining formulas: one with a block of existential quantifiers followed by a formula
with only bounded quantifiers, $ \forall x <y$ and $ \exists x <y$, and the other one with a block of universal quantifiers followed by a formula with only bounded quantifiers. A block of existential
(universal) quantifiers can be replaced by a single existential (universal) quantifier by coding tuples of natural numbers by a single natural number. It
follows from the proof of Hilbert's Tenth Problem that bounded quantifiers can be eliminated from the above formulas, so a
computable set is definable in $\mathcal{N}$ both by an existential and a universal formula. 


\begin{figure}[htb!]
\centering

    \includegraphics[trim = 0.4cm 9.7cm 0.2cm 0.1cm, clip, width=0.44\textwidth]{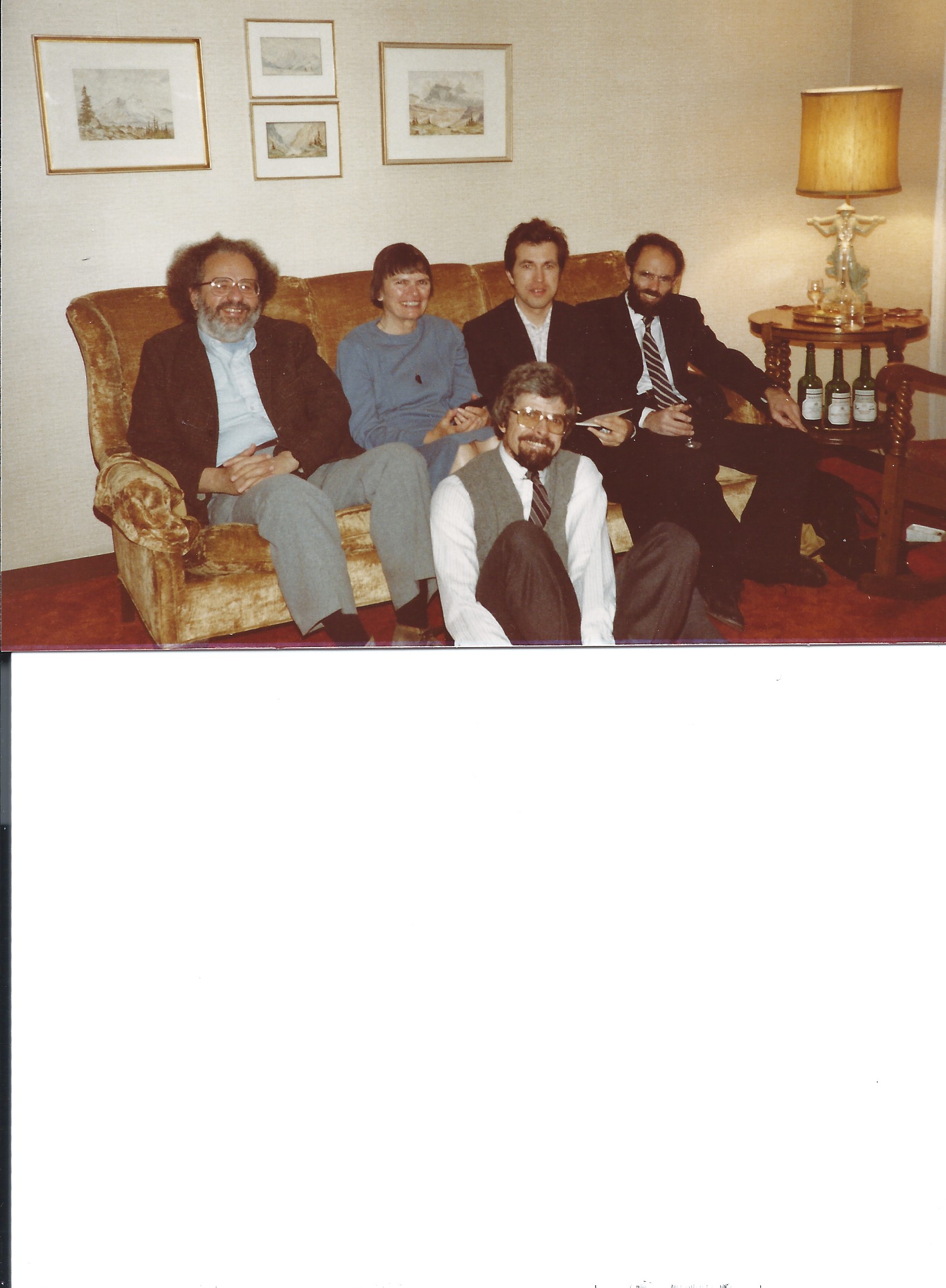}
    
\caption{\label{fig:MartinJuliaYuri}Martin Davis with Julia Robinson, Yuri Matiyasevich, James P. Jones (and Patrick Brown, front) in Calgary, 1982.}
\end{figure}

\subsection*{``Positive aspects of a negative solution''}
This was a part of the title of the 1974 paper by Davis, Matiyasevich and Robinson \cite{DMR76}.  Perhaps this title was a response to some opinions in the mathematical community that the negative answer to Hilbert's question about polynomial equation meant that the subject matter was closed.  Nothing could have been further from the truth. 

Among other things, the authors of the 1974 article explained that in a manner of speaking the negative solution was inevitable in part because a big  part of Mathematics can be encoded into polynomial equations, e.g., Riemann's Hypothesis.  More precisely, Riemann Hypothesis is true if and only if a certain polynomial equation with known coefficients  has no integer solutions!  Such polynomials exist for many other famous problems: Goldbach's conjecture, consistency of ZFC, etc.  Thus existence of an algorithm to solve polynomial equations would resolve unreasonably many open questions in mathematics.

Perhaps the most positive consequence of Davis--Putnam--Robinson--Matiyasevich theorem is that it, together with definability results of Robinson, became a foundation of a  new field: definability and decidability in number theory. This field seeks to understand what is definable and decidable in the language of rings (i.e., the language of polynomial equations) and its various dialects over rings and fields of interest to number theory.  From its inception this field was situated on the boundary of several mathematical areas: number theory, algebraic geometry, model theory and computability theory and led to some interesting interactions between these fields.
\subsection*{The question of $\Q$}

Perhaps the most important question in this new area is the analog of Hilbert's Tenth Problem  for $\Q$.  More precisely, the question is whether there exists an algorithm that can determine whether an arbitrary polynomial equation in several variables with integer coefficients has solutions in $\Q$. One can show that a positive answer to Hilbert's question for $\Z$ implies a positive answer to the question over $\Q$.  However, the reverse implication is not clear.  

One way to show that there is no algorithm to determine whether polynomial equations have solutions over $\Q$ is to construct a Diophantine definition of $\Z$ over $\Q$.  However, there are conjectures by Mazur and others implying that such a definition does not exist.  The question concerning (non)existence of this Diophantine definition is another major problem in the area.

\subsection*{On a personal note AS}
I first spoke to Martin (on the phone) in the Spring of 1983 when I was deciding on a graduate school.  Martin encouraged me to come to NYU.  While there, I was lucky enough to take his class on Hilbert's Tenth Problem.  Martin developed a different method for showing that exponential equations were Diophantine (polynomial) using the Pell equation in place of Fibonacci numbers.  His method and its generalizations to norm equations served me well in many a paper.  

I continued to be in touch with Martin until his death seeking his advice on many issues.  He was mathematically engaged until the very end.  I believe his last talk took place online.  It was the opening  talk for the MSRI semester on Definability, Decidability and Computability in Number Theory mentioned above.
\begin{figure}[!htbp]
    \centering
    \includegraphics[width=0.4\textwidth]{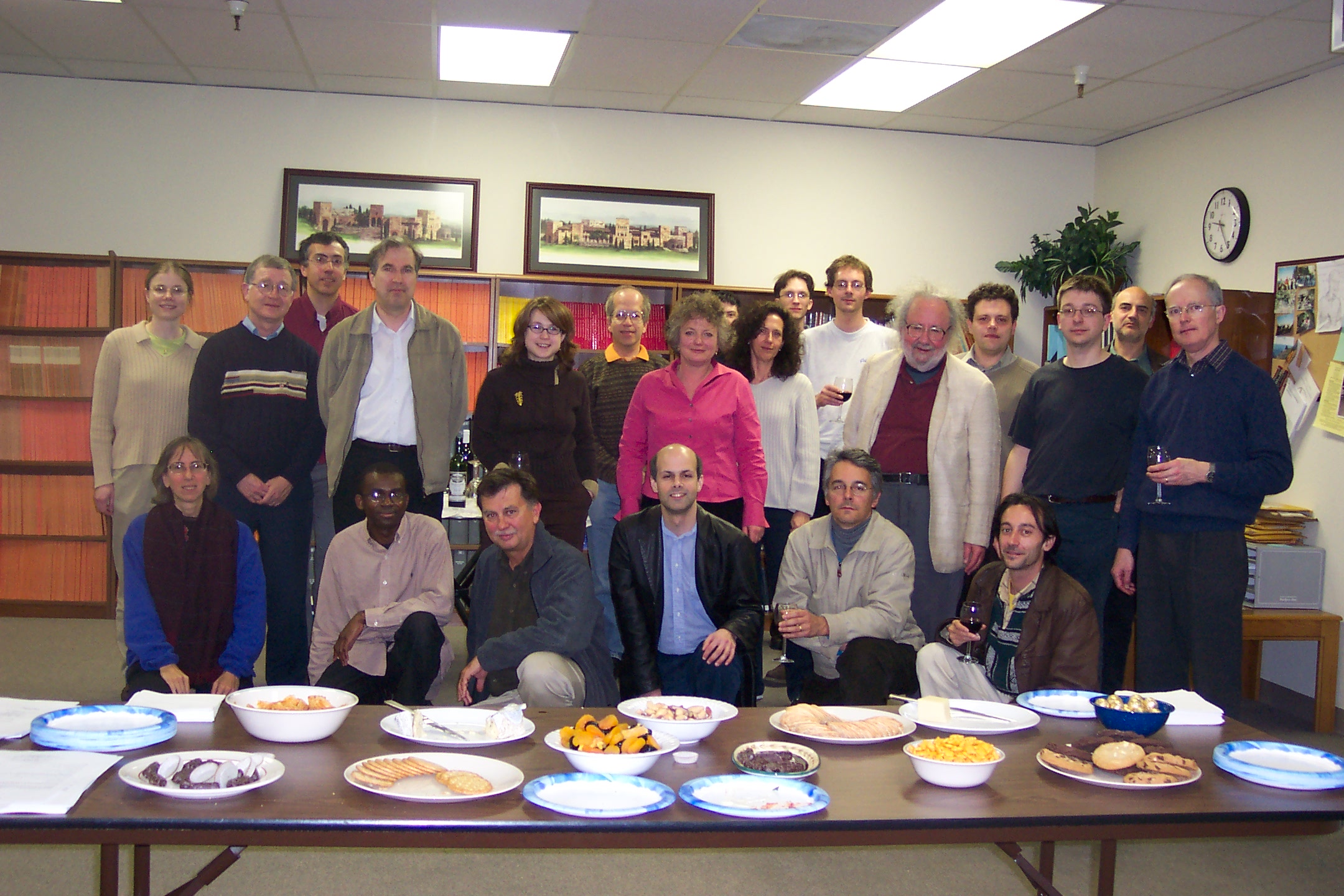}
    \caption{Participants at a meeting at the American Institute of Mathematics on extensions of Hilbert's Tenth Problem, 2005.}
    \label{fig:AIM}
\end{figure}

\section{Automated Reasoning}
\begin{sloppypar}
\noindent    Martin began to program computers in 1951, when he was recruited in a group\HIDE{\footnote{In \cite{Dav16} Martin recalls that it was ``the hot breath of the draft'' that led him to leave his ``Research Instructor'' position at the University of Illinois at Champaign-Urbana. His next position was in the organization, led by Frederick Seitz, named Control Systems Laboratory. Thanks to this move and, subsequently, enjoying a two-year ONR grant at the IAS in Princeton, he managed to avoid being inducted into the army.}} that developed programs for an \textsc{ordvac} machine, in support of the military during the Korean War \footnote{Martin then entered the organization, led by Frederick Seitz, named Control Systems Laboratory. Thanks to this move and, subsequently, enjoying a two-year ONR grant at the IAS in Princeton, he managed to avoid being inducted into the army. Concerning \textsc{ordvac} and \textsc{johhniac}, cf. \cite[p.$\:${162}]{Davis18a}.}. He was assigned the task of writing, in absolute binary machine language, the prototype of a system by which \textsc{ordvac} was to navigate 100 airplanes in real time. After this tumultuous training on concrete programming, that lasted roughly one year,  Martin was confident enough on his skills with computers that he managed to receive funding for a project on implementing Presburger's decision algorithm for integer arithmetic without multiplication. Martin's implementation of that algorithm took place on a \textsc{johnniac} computer available at the Institute for Advanced Study in Princeton in the summer of 1954. The prover could only ascertain very simple statements  \cite{Davis60,Davis01a}, yet its accomplishment marked a milestone in computational logic. The preface of \cite{SiekmannWrightson83a} says that it produced ``what appears to be the first computer generated mathematical proof,'' and this accomplishment qualifies Martin as a trailblazer of the field today known as `automated reasoning.' For some twenty-five years,\footnote{Still much later, around 1990, Martin would again look at the problem of automatic proof discovery. In the paper \cite{Davis-Fechter-91} he coauthored, that presents first-order predicate calculus under a very unusual light, he points out: ``The above very elementary examples only hint at the kinds of proof procedures which our free variable formulation should make possible. But there is reason to believe that they may turn out to be of interest''.}  Martin continued to contribute to that field.

In the late 1950s, the seminal report \cite{DavisPutnam58a} on computational methods for the propositional calculus arose from his collaboration with the distinguished analytic philosopher Hilary Putnam. With him Martin enjoyed discussing ``all day long about everything under the sun'' for various consecutive summers. The Davis-Putnam-Logemann-Loveland procedure \cite{DavisLogemannLoveland62} (DPLL for short), still fundamental in today's architectures of fast Boolean satisfiability solvers, was rooted in that collaboration. Between 1958 and 1960, simultaneously with the Davis-Putnam and DPLL projects, three major projects (led by Gilmore, Dunham-Fridshal-Sward, and Wang, respectively---see \cite{LSS16}) were developing propositional provers. It was Davis and Putnam who set up the overall organization that, after them, would become standard in the automation of quantification theory. They adopted the Conjunctive Normal Form (CNF, a propositional conjunction of disjunctions of affirmed or denied logical variables) in pursuing an unsatisfiability test and embedded their rules for propositional logic into the enhanced proof framework.  Martin and Hilary Putnam viewed a tester able to establish whether or not a given CNF formula is truth-functionally satisfiable as a key component in a general-purpose procedure for quantification theory.\footnote{Let us recall that in quantification theory, unlike propositional logic, the problem of validity is semi-decidable (that is, computably enumerable) without being algorithmically solvable. That is to say: while a systematic search will sooner or later unearth the proof of a theorem, rejecting an unprovable conjecture may turn out impossible.}   This general procedure can then be applied to obtain proofs (by contradiction) in virtually any mathematical domain.  It is surprising that, to these days, the DPLL-procedure constitutes
 a kernel of any efficient CNF-tester.

Any propositional formula can be brought to an equisatisfiable
CNF formula in linear time (cf. \cite{ts70}), hence CNF-satisfiability would become paradigmatic of the whole collection of $\mbox{NP}$-hard problems.\footnote{Cf. \url{https://www.claymath.org/wp-content/uploads/2022/06/pvsnp.pdf}.} 
Martin never took sides on the ``\emph{$P$ vs $NP$}'' problem: he argued that we have no reliable intuition of what an algorithm of, say, complexity $n^{27}$ can do.  He discussed his views on the problem at more length in \cite{Davis20}.

Martin and Hilary Putnam had made it clear that their method would outperform competitors of the time: by exploiting it in a 30-minute \emph{hand computation}, they in fact validated a claim that a theorem-proving program developed by Paul C.\ Gilmore had been unable to validate with a 21-minute run on an IBM~704 machine \cite{DavisPutnam60}.  The improvement was not due to a better handling of quantifiers, but due to an improvement in the propositional part \cite{PPV60}. 
 Later, the promise of these hand computations was realized in computer implementations.  Davis and Putnam's proof procedure for finitely axiomatized theories was implemented by Logemann and Loveland at New York University.  They found and removed a bottleneck in the propositional-level component of the procedure.  Later an implementation in LISP at Bell Labs gave further incremental improvements.\footnote{The mentioned implementation at NYU improved the 30 minute hand computation time to 2 minutes (cf. \cite{DavisLogemannLoveland62}).} 

The Bell Labs implementation of Martin's method \cite{Davis63b,Chinlund-Davis-etAl64,Yarmush76} was named \emph{Linked Conjunct}.  The operating principle required that each logical variable in an unsatisfiable CNF formula would be paired with the same variable with the opposite sign in another conjunct. Not long after, there would be a proliferation of new proof search methods arising from John Alan Robinson's influential resolution principle\footnote{Cf. \url{https://www.programmazionelogica.it/2016/11/john-alan-robinson/}}. It turns out, however, that many of these refinements can be naturally explained from the standpoint of Linked Conjunct \cite{Omo82}.

\begin{figure}[htb!]
    \centering
    \includegraphics[width=0.4\textwidth]{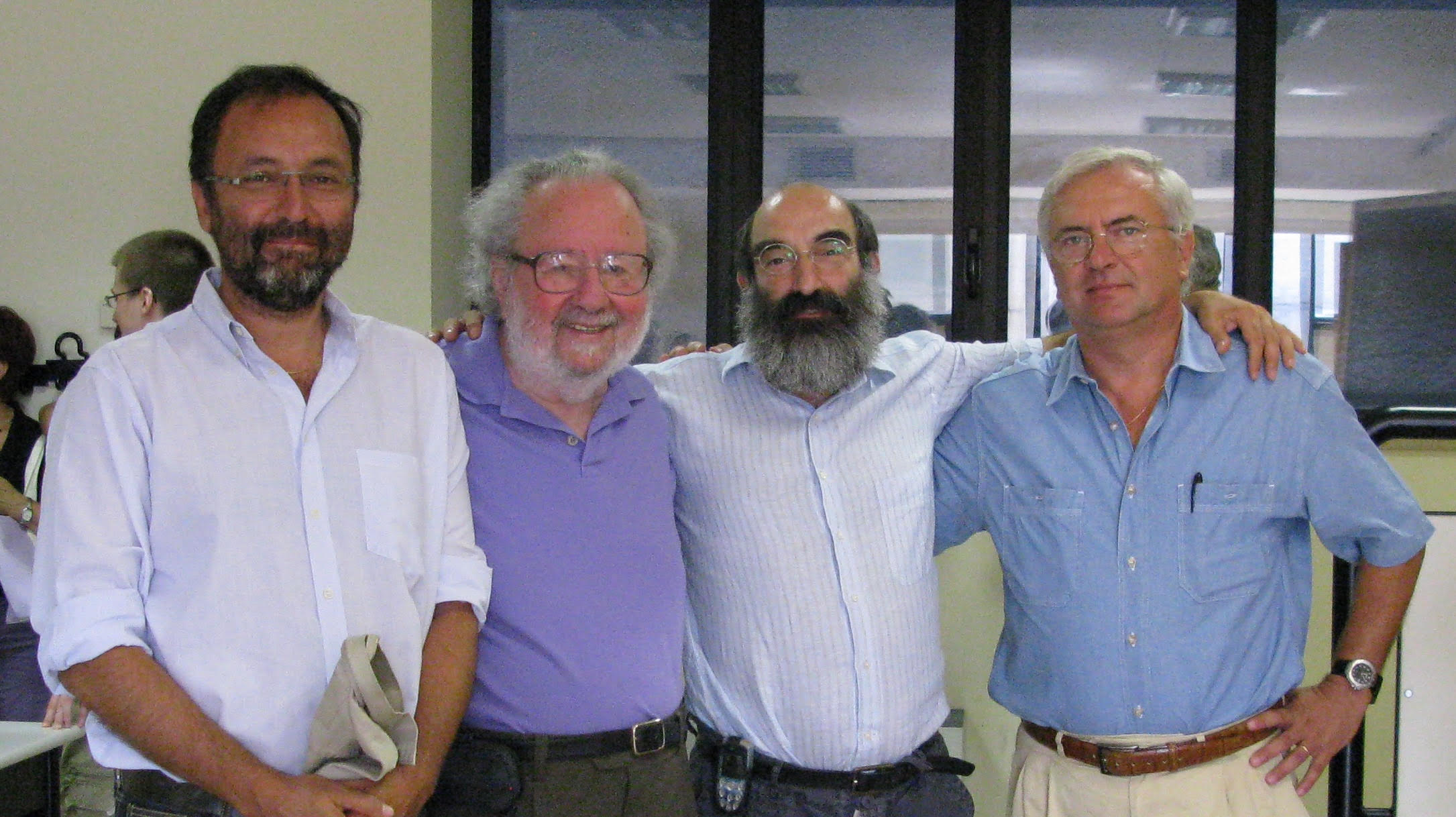}
    \caption{Alberto Policriti, Martin Davis, Eugenio Omodeo, and Franco Parlamento in Trieste in 2009.}
    \label{fig:MartinEugenio}
\end{figure}

The original expectations about stand-alone theorem provers have been somewhat retargeted, over the years, to proof-checking systems that range from highly interactive reasoning assistants to mere proof-script verifiers. Martin also had a role in this change of perspective, as witnessed by his papers \cite{davis-schwartz79} and \cite{Davis81obviouslogical}. The former of these, jointly authored by him and his lifelong friend and colleague ``Jack'' (namely Jacob T. Schwartz, cf. \cite{martinAndEdOnJack,collectionForJack,JackMeetsMarx}), addressed the issue of \emph{metamathematical extensibility} in a full-blown program- and proof-verification technology. Which mechanisms can ensure long-term reliable use 
of a proof checker which undergoes augmentations with new symbols, schemes of notation, and extended rules of inference? 

As for \cite{Davis81obviouslogical}, it stemmed from Martin's experimentation with Richard Weyhrauch's FOL proof checker developed at John McCarthy's Artificial Intelligence Laboratory.\footnote{Martin would later cooperate, with his \cite{Davis80b} conceived in the same stimulating environment at Stanford University, to the launch of non-monotonic logic formalisms.} Martin recounts in \cite{Dav16}: ``I found it neat to be able to sit at a keyboard and actually develop a complete formal proof, but I was irritated by the need to pass through many painstaking tiny steps to justify inferences that were quite obvious'', and then adds: ``Using the LISP source code for the linked-conjunct theorem prover\ldots, a Stanford undergraduate successfully implemented an `obvious' facility as an add-on to FOL''.

\subsection*{On a personal note EO} In 1975 Martin offered a summer course on computability in the pretty Italian town of Perugia. The dozen students in his class were initially amazed at the discrepancy between Martin as an unpretentious, easygoing person and his reputation as a distinguished scholar. Admiration quickly prevailed over astonishment when Martin began his lectures: for the entire one-month duration of his course, concepts remained clear, precise, and accessible. Even when he reached his cherished advanced topic: Hilbert's Tenth Problem.  

I owe to having been a student in Martin's class in Perugia the fact that I could complete my academic formation at NYU, where Martin was my advisor for a Master's and next for a Ph.D. degree in Computer Science. A stream of Italian students and researchers (three had been my classmates in Perugia; Alberto Policriti and others belong to a successive generation) would, like me, reach Martin overseas in the following decade. This testifies to the influence that Martin's crystal-clear lectures, and the subtlety with which he addressed issues of philosophical relevance and depth, used to exert on his audience --- in Italy much as elsewhere. 

\subsection*{On a personal note AP}

On the evening of a beautiful day of the fall of 1990 I was invited by Martin at his place in the upper-west side in New York, to a ``party for two Yuri's''. The two Yuri's were Yuri Gurevich and Yuri Matiyasevich, both in town and hosted by Martin. For Yuri Matiyasevich that was the first visit to the United States. 

The number and (even more) the names of the people invited to the party were --- especially for a young Italian Ph.D. student --- rather overwhelming, and I was definitely scared when I  was greeted amicably by Virginia upon my arrival. However, every tension promptly dissolved when Martin introduced Yuri Matiyasevich to the audience. He started recalling what remained to be proved after his first reduction of Hilbert's Tenth Problem and how, at the time, he conjectured that the last pending issue  (a  number theoretic hypothesis raised by Julia Robinson around 1950) would have certainly be solved by a \emph{\ldots clever young Russian \ldots} in the  near future. 

Yuri was then introduced to everybody as the (constructive and positive) proof of his conjecture.


\begin{figure}[!htbp]
    \centering
    \includegraphics[width=0.3\textwidth]{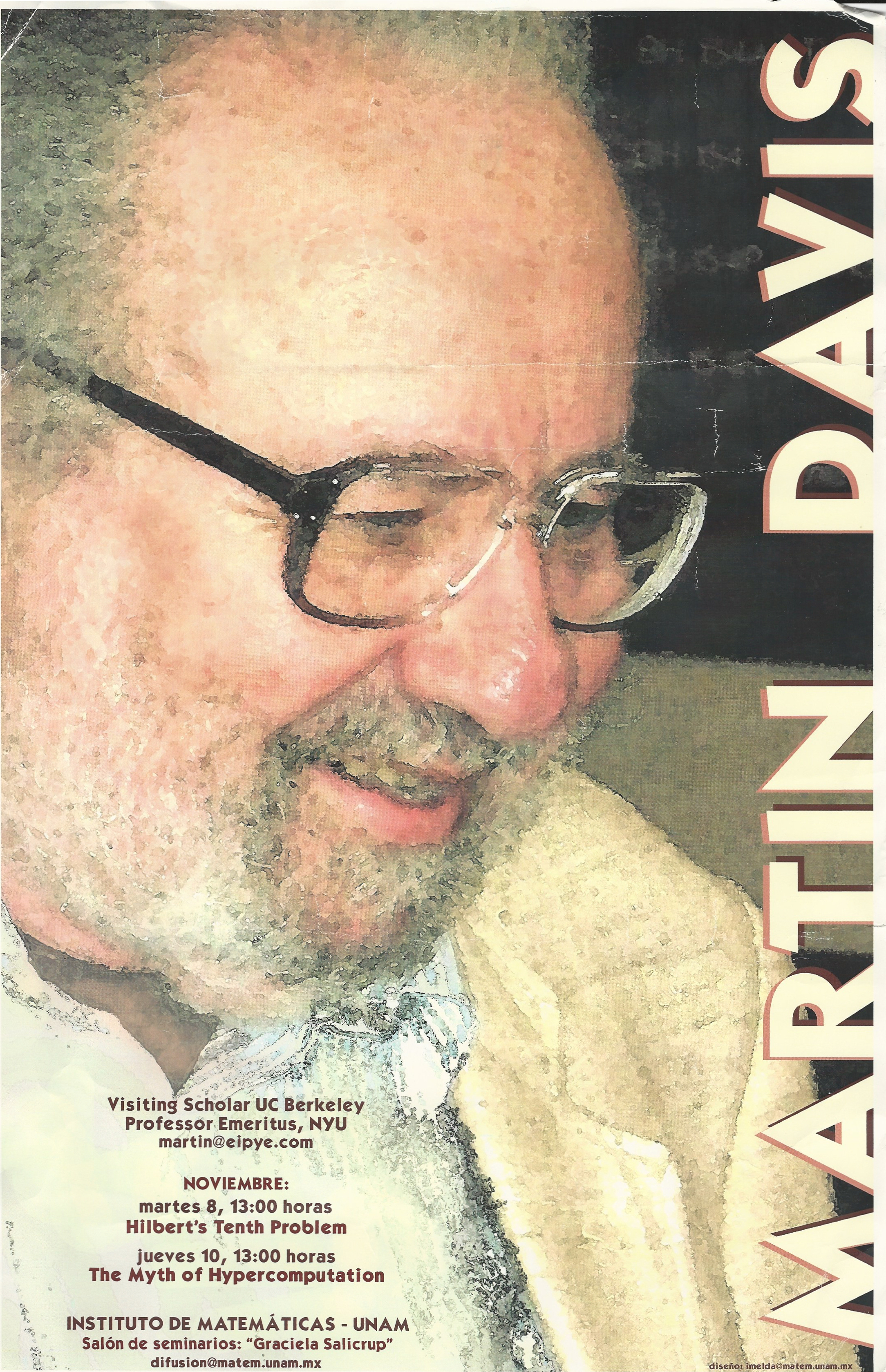}
    \caption{Announcement of a conference by Martin in Mexico. (Courtesy of Laura Elena Morales Guerrero)}
    \label{fig:MartinInMexico}
\end{figure}

\end{sloppypar}
\section{Conclusions}
Recollections of contributions by Martin to computability theory, Hilbert's Tenth Problem, and Automated Reasoning have been scattered all over the preceding text; many more could be cited: e.g.,
in \cite{Davis72},  Martin stretches the algorithmic unsolvability of Hilbert's Tenth Problem into this result:
For each proper nonempty subset $A$ of $\mathbb{N} \cup \{\aleph_0\}$, 
no algorithm can establish, given any polynomial $p$ with integer coefficients, whether the number of distinct positive integer solutions to the equation $p = 0$ belongs to $A$.

Beyond the proofs of specific theorems, Martin's scientific legacy included a broader contribution in the promotion of formal methods and theoretical computer science.  He significantly contributed to the recognition of computability theory as an autonomous branch of mathematics. Martin developed a program in logic, and formed a logic group, first at Yeshiva University (being able to involve in it stars such as Raymond Smullyan) and then at the Bronx campus of the Courant Institute (NYU). In the 1960s, he  devoted a good deal of time and energy in preparing an anthology of fundamental articles by G\"odel, Turing, Post, Kleene, and Rosser \cite{Davis65}.  His involvement with symbolic logic originated in the 1940s from his passionate interest in the foundations of real analysis, which also led him to author a classic on nonstandard analysis in the 1970s and to serve for decades as the moderator of FOM, an automated e-mail list for discussing foundations of mathematics (see \url{https://cs.nyu.edu/mailman/listinfo/fom}). 

Over the years, Martin lectured in several countries (to cite a few: Italy, Japan, India, England, Russia, and Mexico), and his lectures have---along with his publications---exerted a wide influence. The centennial of Frege’s \textit{Begriffsschrift}, Martin reports, ``fundamentally changed the direction of my work'' \cite{Dav16}: Being invited to place some contemporary trends in a proper historical context, he finds ``trying to trace the path from ideas and concepts developed by logicians \dots to their embodiment in software and hardware \dots endlessly fascinating''  \cite{Dav16}.

\begin{figure}[!htbp]
    \centering
    \includegraphics[width=0.4\textwidth]{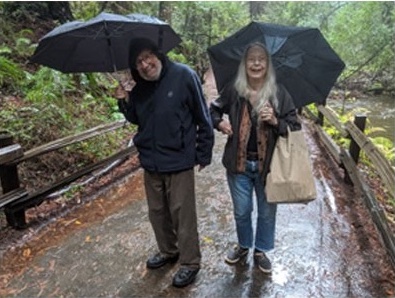}
    \caption{Martin and his wife Virginia in 2019. (Courtesy of Nina and Yuri Matiyasevich)}
    \label{fig:MartinVirginia2019}
\end{figure}

\end{sloppypar}

\bibliographystyle{plain}\footnotesize
\bibliography{MDbibliography}

\end{document}